\pdfoutput=1

\documentclass[12pt]{article}%
\usepackage{amssymb,amsmath,xcolor,graphicx,xspace,colortbl,ragged2e,rotating}
\usepackage{amsfonts}
\usepackage{amsmath}
\usepackage{amssymb}
\usepackage{amsthm}
\usepackage{boxedminipage}
\usepackage{color}
\usepackage{graphics}
\usepackage{graphicx}
\usepackage{ragged2e}
\usepackage{url}
\usepackage{xcolor}
\usepackage{hyperref}%
\providecommand{\U}[1]{\protect\rule{.1in}{.1in}}
\graphicspath{{how-did-fermat-discover-his-theorem-24_graphics/}{how-did-fermat-discover-his-theorem-24_tcache/}{how-did-fermat-discover-his-theorem-24_gcache/}}
\DeclareGraphicsExtensions{.pdf,.eps,.ps,.png,.jpg,.jpeg}
\providecommand{\U}[1]{\protect\rule{.1in}{.1in}} \providecommand{\U}[1]{\protect\rule{.1in}{.1in}}

\newtheorem*{exercise-unnumbered}{Exercise}
\newtheorem*{lemma-unnumbered}{Lemma}
\newtheorem*{corollary-unnumbered}{Corollary}

\newtheorem*{GFTP}{Generalization of Fermat's third proposition}
\newtheorem*{DC}{Divisibility conjecture}
\newtheorem*{FBT}{Fermat's broader theorem}
\newtheorem*{FT}{Fermat's theorem}
\newtheorem*{FFP}{Fermat's first proposition}
\newtheorem*{FSP}{Fermat's second proposition}
\newtheorem*{FTP}{Fermat's third proposition}
\newcommand{\sep}{\vspace{-3pt} \begin{center} {\mathversion{normal}
$\infty \mspace{-5.5mu} \infty \mspace{-5.5mu} \infty \mspace{-5.5mu} \infty \mspace{-5.5mu}
\infty \mspace{-5.5mu} \infty \mspace{-5.5mu} \infty \mspace{-5.5mu} \infty$}
\end{center} \vspace{-3pt}}
\begin{document}

\title{How did Fermat discover his theorem? }
\author{David Pengelley }
\maketitle

What is perhaps both the first and the most important surprising property ever
discovered about prime numbers? We maintain it is Pierre de Fermat's theorem
on power residues modulo a prime, which he asserted in 1640. \cite[ch.
I,\S 2,3; ch. II,\S 4]{weil}

\begin{FT}
If $p$ is prime, and $a$ is not divisible by $p$, then $a^{p-1}-1$ is
divisible by $p$. (In the language of power residues, $a^{p-1}$ has residue
$1$ modulo $p$).
\end{FT}

We ask here: How could Fermat possibly have discovered something so surprising?%

\begin{center}
\includegraphics[
height=2.4278in,
width=1.7141in
]%
{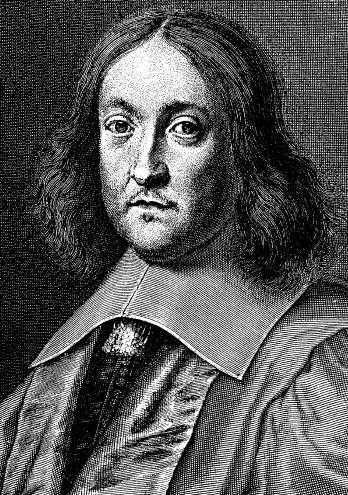}%
\\
Pierre de Fermat (1607--1665)
\end{center}

Two millennia prior to Fermat, two other key features of primes had emerged,
appearing in Euclid's \emph{Elements}, based on the definition of a prime as
multiplicatively indecomposable. One is their infinitude in number, but which
is not that unexpected. The other is \emph{Euclid's lemma}, that if a prime
divides a product then it divides one of the factors. Euclid's lemma is
critical not only to number theory, but to all of mathematics, since it
ensures the uniqueness of prime factorization in the natural numbers, which
mathematicians use frequently whether they consciously realize it or not. But
even though proving Euclid's lemma is actually much subtler than Euclid
incorrectly thought \cite{pengelley-richman,pengelley-fractions}%
\cite[\S 4.4]{SG-book}, and the result is a wonderful and incredibly fortunate
feature of the natural numbers, it is not really all that unexpected either.

Fermat's theorem, on the other hand, is totally unexpected and truly
astonishing. There is also little in number theory that is not dependent on it
or intertwined with it. Its remarkable importance even today is amply
demonstrated by the fact that now, almost four centuries later, Fermat's
theorem provides the mathematical foundation for the RSA cryptosystem, which
is still central to society's communications security even after several
decades serving as its heart. \cite[\S 7.2]{SG-book}\cite{rsa,rsa-wikipedia}

We will investigate here how Fermat discovered his theorem, largely because it
is so unobvious. Fermat's pursuit at the time is well known, but imagining
exactly how this led him to guess his theorem is our passion here. We will
also see some lessons about the process and dissemination of discovery.

Along the way we will provide a few exercises, giving a tiny glimpse of how
one might use this development and analysis in creating a teaching project for
a number theory course. See \cite[appendix A]{SG-book} for a detailed more
advanced set of inquiry-based activities. 

\subsection*{Perfect Numbers}

Fermat was as captivated as others of his day with perfect numbers, a
classical Greek concept that may seem a mere curiosity today, but which
spawned Fermat's theorem as a monumental consequence. \cite[pp. 81,
497--8]{katz}\cite[ch. VI,\S II]{mahoney}\cite[ch. II,\S 4]{weil}

A \emph{perfect number} is one that equals the sum of its proper divisors (its
\emph{aliquot parts}), e.g., $6=3+2+1$ or $28=14+7+4+2+1$. Already in Euclid's
proposition IX.36 \cite{euclid} we learn how to construct a family of perfect numbers:%

\sep

\textsf{If as many numbers as we please beginning from a unit are set out
continuously in double proportion until the sum of all becomes prime, and if
the sum multiplied into the last makes some number, then the product is
perfect.}%

\sep

You will not find it hard to rephrase this more modernly and to confirm
Euclid's claim: If $1+2+\cdots+2^{n-1}=2^{n}-1$ is prime, then $\left(
2^{n}-1\right)  2^{n-1}$ is perfect.

So there is a pairing between Euclid's perfect numbers and prime numbers of
the form $2^{n}-1$. Before pursuing finding such primes, as did Fermat and his
contemporaries, we remark that Leonhard Euler later showed that Euclid's
family captures all the even perfect numbers \cite[v. 1, ch. I]{dickson}%
\cite[E798]{eulerarchive}\cite[pp. 81, 497--8]{katz}\cite[p. 182]{SG-book}.
But even today we do not know if there are infinitely many numbers in Euclid's
family, because we do not know if there are infinitely many primes $2^{n}-1$
\cite{gimps}. Nor do we know whether there are any odd perfect numbers at all.

\begin{exercise-unnumbered}
Find the third even perfect number.
\end{exercise-unnumbered}

Fermat and his correspondents were also interested in other similar phenomena
centering on the sum of the divisors of a number. This led them to develop
methods to analyze more generally when numbers of the form $a^{n}-1$ or
$a^{n}+1$ are prime, or what their prime divisors could be, a further
fascinating journey. \cite[ch. VI,\S II]{mahoney}\cite[appendix A]%
{SG-book}\cite[ch. II,\S 4]{weil}

\subsection*{Mersenne Numbers and Fermat's Broader Theorem}

Thus Fermat and his correspondents, with Marin Mersenne as frequent gobetween
\cite[ch. VI,\S II]{mahoney}\cite[ch. II,\S 4]{weil}, were very interested in
understanding the potential primality or divisors of any number of the form
$M_{n}=2^{n}-1$, which today is called a \emph{Mersenne number}, or a
\emph{Mersenne prime} if it is prime\footnote{The largest primes we know today
are Mersenne primes.}. We will, with Fermat, call $n$ the \emph{exponent}\ of
$M_{n}$. It was this primality investigation that somehow led Fermat to his
theorem, and our goal is to glean from his letters how Fermat's path of
discovery may well have proceeded.

Most who quote Fermat on his theorem do so from his letter of 18 October 1640
to Bernard Frenicle de Bessy, one of his keenest mathematical correspondents,
so let us examine that first. \cite[pp. 206--212]{fermatoeuvres}\cite[ch.
VI,\S II]{mahoney}\cite[ch. II,\S 4]{weil} He called it the \textquotedblleft%
\textsf{fundamental proposition} \textsf{ of aliquot parts}\textquotedblright%
\ and \textquotedblleft\textsf{the fundament on which I\ support the
demonstrations of everything concerning geometric progressions}%
\textquotedblright. \cite[pp. 208--209]{fermatoeuvres} Clearly Fermat judged
this proposition as foundationally important, and key to studying perfect numbers.

Fermat presented his proposition to Frenicle thus:%

\begin{center}
\includegraphics[
height=1.1028in,
width=3.9054in
]%
{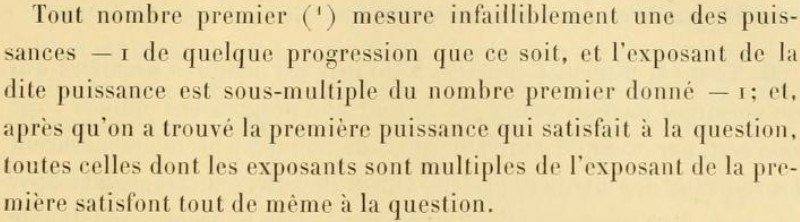}%
\\
Fermat's (broader) theorem, 18 October 1640 \cite[p. 209]{fermatoeuvres}%
\end{center}
%

\sep

\textsf{Without fail, every prime number measures one of the powers }$-1$
\textsf{of any progression whatever, and the exponent of the said power is a
submultiple of the given prime number }$-1$.\textsf{ Also, after one has found
the first power that satisfies the problem, all those of which the exponents
are multiples of the exponent of the first will similarly satisfy the problem.
}%

\sep

Let's interpret this a little bit, place it in modern terminology, and see if
it sheds any light on how Fermat discovered his theorem.

By \textquotedblleft measures\textquotedblright\ Fermat means
\textquotedblleft divides\textquotedblright, and by \textquotedblleft
submultiple\textquotedblright\ he means \textquotedblleft
divisor\textquotedblright. By \textquotedblleft progression\textquotedblright%
\ he means a geometric progression in which the terms are the powers of a
fixed number, called the base. And when Fermat writes \textquotedblleft%
\textsf{every prime number measures one of the powers } $-1$ \textsf{ of any
progression whatever, and the exponent of the said power is a submultiple of
the given prime number } $-1$ \textquotedblright, the \textquotedblleft%
\textsf{said power}\textquotedblright\ refers to the first such power. Also,
Fermat didn't really mean \textquotedblleft\textsf{any progression}%
\textquotedblright, since obviously it isn't true for the prime $3$ with the
progression consisting of powers of $6$. He had in mind that the base for the
progression is not divisible by the chosen prime number. Finally, Fermat meant
not only that \textquotedblleft\textsf{all those of which the exponents are
multiples of the exponent of the first will similarly satisfy the
problem}\textquotedblright\footnote{By \textquotedblleft\textsf{satisfy the
problem\textquotedblright}, Fermat means that the \textquotedblleft%
\textsf{prime number measures one of the powers }$-1$\textquotedblright.}, but
that these are the only exponents that satisfy the problem; we will discuss
this subtlety shortly, and see soon after why it is important. Putting this
all together, we can state what Fermat claimed thus.

\begin{FBT}
[October 1640, according to Fermat]If $p$ is prime, and $a$ is not divisible
by $p$, then there is a least positive exponent $k$ for which $a^{k}-1$ is
divisible by $p$. And for any $n$, $a^{n}-1$ is divisible by $p$ precisely
when $n$ is a multiple of $k$. Moreover, $k$ is a divisor of $p-1$.
\end{FBT}

Henceforth, to keep things interesting, correct, and clear, we will assume
that $p>2$, that $a>1$, and that the exponent $n$ in numbers of the form
$a^{n}-1$ is positive.

We call Fermat's statement here a \textquotedblleft
broadening\textquotedblright\ because it claims more than what we label as
Fermat's theorem today. At first sight the broadening might appear
considerably stronger, but it is easily seen to be equivalent to Fermat's
theorem as we know it today (see below). The interesting hint here for our
question \textquotedblleft How did he discover it?\textquotedblright\ is that
Fermat is discussing not only the $p$-divisibility of $a^{p-1}-1$ for $p$
prime, but analyzing the potential $p$-divisibility of $a^{n}-1$ for all $n$.
In other words, for $a=2$, this is about all possible prime divisors of all
Mersenne numbers, which matches with what we know he had been studying for
pursuing perfect numbers.

Fermat's claim here is twofold. First, a given [odd] prime $p$ will divide
precisely all the Mersenne numbers (or more generally all $a^{n}-1$ for fixed
$a$ not divisible by $p$) whose exponents $n$ are multiples of a smallest
positive exponent $k$, i.e., $\left\{  n=jk:j=1,2,\ldots\right\}  $. We will
call such an arithmetic progression, in which the first term is also the
increment, a \emph{simple arithmetic progression}. Second, Fermat claims that
this progression will include $p-1$.

We will see shortly why the broadened version Fermat stated here is highly
relevant to divining how he discovered his result.

First, though, since we encompassed in our interpretation above of Fermat that
the multiples $n$ of the claimed $k$ are the only exponents for which
$a^{n}-1$ is divisible by $p$, and since Fermat didn't literally say this,
let's make sure that this would have been straightforward for him to include
implicitly\footnote{We could use congruence notation, or residues, or group
theory, throughout to ease our way, but this would be anticipating by a
century and a half or more, and we would miss experiencing the authenticity of
Fermat's path. So we refrain.}.

\begin{lemma-unnumbered}
a) If $k$ divides $n$, then if $p$ divides $a^{k}-1$, also $p$ will divide
$a^{n}-1$.

b) Essentially conversely, if $p$ divides $a^{n}-1$, and if $k$ is the least
exponent for which $p$ divides $a^{k}-1$, then $k$ divides $n$.
\end{lemma-unnumbered}

%

\begin{proof}%

Fermat would have been familiar with the standard algebraic identity involving
a geometric sum:
\begin{equation}
a^{st}-1=\left(  a^{t}-1\right)  (a^{(s-1)t}+a^{(s-2)t}+\cdots+1).
\label{geometric-sum-identity}%
\end{equation}

The consequent divisibility relation $(a^{t}-1)\mid(a^{st}-1)$ immediately
confirms part a), his claim that \textsf{\textquotedblleft all those of which
the exponents are multiples of the exponent of the first will similarly
satisfy the problem\textquotedblright} (i.e., be divisible by $p$).

To see part b), that these are the only exponents that satisfy the problem, we
consider any $a^{n}-1$ divisible by $p$, and write $n=jk+r$ with $0\leq r<k$.
Our goal is to show that $r=0$.

We write $a^{n}-1=a^{jk}\cdot a^{r}-1=(a^{jk}-1)a^{r}+(a^{r}-1)$. Now since
$p$ divides $a^{k}-1$, and thus divides also $a^{jk}-1$, as well as $a^{n}-1$
by hypothesis, we deduce that $p$ divides $a^{r}-1$. But $k$ is the least
positive exponent with this property, so we must have $r=0$, as claimed.%

\end{proof}%

\begin{corollary-unnumbered}
If $p$ divides some Mersenne number (or more generally some $a^{n}-1$ for
fixed $a$ not divisible by $p$), then $p$ will divide precisely all such
numbers whose exponents $n$ are multiples of a smallest positive exponent $k$,
i.e., $\left\{  n=jk:j=1,2,\ldots\right\}  $.
\end{corollary-unnumbered}

N.B: The Lemma and Corollary do not require that $p$ be prime, and we will use
this later. Primality is required only for Fermat's further claim that $k$
exists and divides $p-1$.

In Fermat's October 1640 letter stating his broader theorem, there is no
indication of how he arrived at his claims; he adds only \textquotedblleft%
\textsf{I would send you the proof, if I did not fear being too long}%
\textquotedblright! \cite[p. 209]{fermatoeuvres} So where else can we look to
see how Fermat discovered this amazingly unexpected result, and to what end?

Earlier than his October letter to Frenicle, Fermat, in his quest to answer
questions about perfect numbers, had discovered an amazing tool to help him
test Mersenne numbers for divisibility by primes. It allowed him to go way
beyond what can be done by hand with just brute force and the sieve of
Eratosthenes. Let us proceed to see a challenge about perfect numbers, and
Fermat's response.

\subsection*{Frenicle's Challenge}

It was typical for the number theory pioneers to issue challenges to each
other in their correspondence, and not to divulge all their methods or proofs.
In March of 1640, Frenicle had written via Mersenne to Fermat, asking for a
perfect number that has at least 20 digits \cite[pp. 182--5]{fermatoeuvres}%
\cite{fletcherFT,fletcher1640corresp}\cite[p. 293]{mahoney}\cite[p. 54]{weil}!%

\begin{center}
\includegraphics[
trim=0.000000in 0.000000in -0.002165in 0.000000in,
height=0.3782in,
width=3.85in
]%
{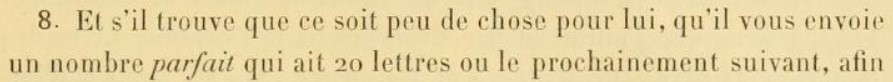}%
\\
Frenicle's challenge to Fermat, March 1640 \cite[pp. 182--5]{fermatoeuvres}%
\end{center}

It is not surprising that Frenicle might issue such a huge challenge, since in
1638 Fermat had boasted that he had a \textsf{\textquotedblleft
method\textquotedblright\ }\cite[pp. 176--7]{fermatoeuvres}\cite[p.
290]{mahoney}\textsf{, }secret of course, for dealing with \textquotedblleft%
\textsf{all questions about aliquot parts}\textquotedblright\ (an aliquot part
meant a proper divisor, and aliquot parts meant their sum, so this included
the problem of finding perfect numbers). Fermat did add, though, that
\textquotedblleft\textsf{the length of the calculations discourages me, as
well as the pursuit of prime numbers, to which all these questions
reduce.}\textquotedblright\ 

\begin{center}%
\begin{center}
\includegraphics[
trim=0.000000in 0.000000in -0.000718in 0.000000in,
height=0.5501in,
width=3.9029in
]%
{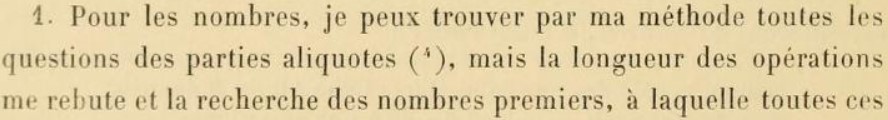}%
\\
Fermat's boast, 26 December 1638 \cite[pp. 176--7]{fermatoeuvres}%
\end{center}

\end{center}

Let's see perhaps why Frenicle challenged Fermat with exactly a 20 digit
minimum, and what Fermat therefore faced in this challenge.

One can immediately restrict, when looking for Mersenne primes, to Mersenne
numbers having prime exponents, since the identity
(\ref{geometric-sum-identity}) above, telling us that $a^{st}-1$ is divisible
by $a^{t}-1$, ensures that composite exponents produce composite Mersenne
numbers. However, a prime exponent is not by itself sufficient to produce a
Mersenne prime: The first exception is $M_{11}=$ $2047=23\cdot89$. Let us call
a Mersenne number $M_{p}$ that is not prime, even though it has a prime
exponent $p$, a \emph{Mersenne prime imposter}. So far we know just one
imposter, $M_{11}\text{,}$ but shortly we will see more, and you should be on
the lookout for patterns.

\begin{exercise-unnumbered}
Observe something interesting about the prime factors of $M_{11}\text{.}$
\end{exercise-unnumbered}

The next few possible Mersenne primes, then, are $M_{13},M_{17},M_{19}%
,M_{23},M_{29},M_{31},M_{37}$. For some of these, but certainly not all, it
was known prior to Fermat whether they were prime or not \cite[v. 1, ch.
I]{dickson}. The perfect number that would correspond to $M_{31}$ (if it were
prime) is $2^{30}M_{31}=2305\,843\,008\,139\,952\,128$, which has $19$ digits.
So Frenicle chose his minimum $20$ digit challenge in order to force Fermat to
look beyond $M_{31}$ to $M_{37}$, whose possible perfect number is
$2^{36}M_{37}=9444\,732\,965\,670\,570\,950\,656$ with $22$ digits.
\cite{fletcherFT,fletcher1640corresp}

\subsection*{Fermat's Shortcuts, Three Very Beautiful Propositions, and a
Great Light}

Thus the first question Fermat faced was whether $M_{37}=2^{37}-1$
$=137\,438\,953\,471$ is prime. But how could he test this? The only known
method would be by brute force division by all primes up to its square root,
$370727.6\ldots$. So he would have to know all the primes up to $370727$, and
then determine divisibility by each of them (there are roughly $30000$ primes
in this range). This would have been way beyond reasonable even for Fermat.
But through his work studying patterns in factorizations of Mersenne numbers,
Fermat had made some astonishing discoveries that could make his task
substantially easier. \cite[ch. VI,\S II]{mahoney}\cite[appendix A]%
{SG-book}\cite[ch. II,\S 4]{weil} In April or May of 1640
\cite{fletcherFT,fletcher1640corresp} he wrote via Mersenne that he had
discovered \textsf{\textquotedblleft several shortcuts\textquotedblright}
(\cite[p. 194]{fermatoeuvres}) for finding perfect numbers, and in June he
revealed them \cite[pp. 194--9]{fermatoeuvres}. (The correspondence between
Fermat and Frenicle, including through Mersenne, is filled with twists, turns,
lacunae, and other unknowns, all fascinatingly investigated, elaborated, and
reconstructed in \cite{fletcherFT,fletcher1640corresp}.) Fermat wrote
\textquotedblleft\textsf{Here are three very beautiful propositions that I
found and proved not without difficulty: I}\textsf{ }\ \textsf{may call them
the foundations of the discovery of perfect numbers. ... From these shortcuts
I already see the birth of a large number of others and this appears to me as
a great light.\textquotedblright}

\begin{center}%
\begin{center}
\includegraphics[
height=0.582in,
width=3.9146in
]%
{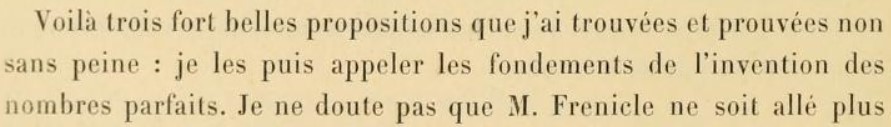}%
\\
Fermat:\ \textquotedblleft\textsf{Three very beautiful
propositions\textquotedblright}, June 1640 \cite[pp. 194--9]{fermatoeuvres}%
\end{center}
%

\begin{center}
\includegraphics[
height=0.3623in,
width=3.9096in
]%
{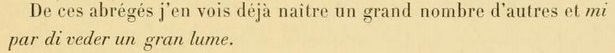}%
\\
Fermat:\ \textquotedblleft\textsf{A great light\textquotedblright}, June 1640
\cite[pp. 194--9]{fermatoeuvres}%
\end{center}

\end{center}

Fermat was clearly extremely pleased with these. We'll write them in modern
terminology, compare them with his October 1640 claim, see how they enabled
him to address Frenicle's challenge, and simultaneously use them to try to
follow Fermat's path in discovering his theorem.

\subsection*{Divisibility Among Mersenne Numbers}

\begin{FFP}
[June 1640]If $n$ is composite, so is $M_{n}$.\footnote{Fermat refers here to
the Mersenne number for a particular exponent as a \textquotedblleft%
\textsf{nombre radical}\textquotedblright, i.e., a \textquotedblleft root
number\textquotedblright; in other words, the type of number he is interested
in for creating perfect numbers.}
\end{FFP}

\begin{center}%
\begin{center}
\includegraphics[
height=0.5946in,
width=3.8928in
]%
{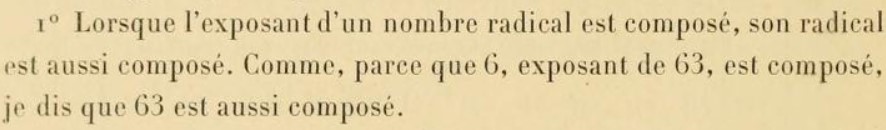}%
\\
Fermat's first proposition, June 1640 \cite[pp. 194--9]{fermatoeuvres}%
\end{center}

\end{center}

This is not news to us, using the identity (\ref{geometric-sum-identity})
above, but perhaps Fermat stated it here because it was extremely helpful in
ways he does not divulge; remember, he and his correspondents were always
challenging one another with results without revealing all their methods and
thinking. Beyond showing that only prime exponents could yield Mersenne
primes, the identity (\ref{geometric-sum-identity}) produces lots of explicit
divisibility relationships between Mersenne numbers, which Fermat could study
for patterns. Specifically, the identity tells us that if $d$ divides $n$,
then $M_{d}$ divides $M_{n}$. This gives great aid in factoring larger
Mersenne numbers using smaller ones. But it also shows that simple arithmetic
progressions of exponents always propagate divisors of Mersenne numbers. This
confirms what Fermat claimed in his \textsf{\textquotedblleft fundamental
proposition\textquotedblright} in October, that the Mersenne numbers divisible
by a given [odd] prime $p$ include all those with exponents in some simple
arithmetic progression. The only caveat to this would be if there were no
Mersenne numbers divisible by $p$. But let's now look at what Fermat says
next, which will resolve the caveat.

\subsection*{Fermat's Theorem in Alternative Form}

\begin{FSP}
[June 1640]If $p$ is [an odd] prime, then $M_{p}-1$ is divisible by $2p$.
\end{FSP}

\begin{center}%
\begin{center}
\includegraphics[
height=0.3849in,
width=3.9146in
]%
{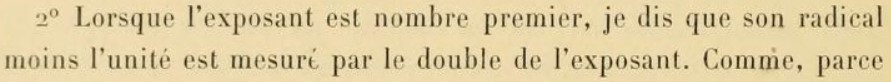}%
\\
Fermat's second proposition, June 1640 \cite[pp. 194--9]{fermatoeuvres}%
\end{center}

\end{center}

This proposition is perplexing at first sight. It is hard to see how it would
be useful for Fermat's goal, to understand the primality or factors of $M_{p}$
for $p$ an odd prime, since it didn't tell him anything directly about factors
of $M_{p}$, even though it is phrased in terms of $M_{p}$. Why would he even
study the factors of $M_{p}-1$ in order to claim something about the possible
factors of $M_{p}$? All the proposition tells us directly about $M_{p}$ is
that it is of the form $2jp+1$. Since numbers of this form are closed under
multiplication, this might recommend proper factors of $M_{p}$ having the same
form, but this is merely suggestive. Why then did Fermat emphasize the second
proposition, and how did he observe it?

The Mersenne numbers satisfy the recursion $M_{n}-1=2M_{n-1}$, so Fermat's
second proposition is equivalent to saying that $M_{p-1}$ is divisible by $p$
(which resolves our caveat above), and in this form it is precisely what we
name as Fermat's theorem today for the base $a=2$! What Fermat stated in his
second proposition then is essentially what we today call the alternative form
of his theorem, that $2^{p}-2$ is divisible by $p$. He had presumably
discovered that $M_{p-1}$ is divisible by $p$, but stated it here in the
alternative form in order to connect it somehow to $M_{p}$. While this
proposition tells us he had discovered Fermat's theorem, it doesn't reveal
anything about how, and the application he makes of the theorem here to
$M_{p}$ does not seem hugely useful for his goal. But we will see that it
presages what he says next.

Fermat's third proposition will be a blockbuster, and it will emerge from our
investigation of how he may have discovered his theorem, which he clearly knew
in June for the base $a=2$. In fact we have seen how the first and second
propositions together provide Fermat's broader theorem for $a=2$: The second
proposition is essentially Fermat's theorem, and the first proposition yields
the broadening to a simple arithmetic progression of exponents.

In what follows, we stick with the base $a=2$, even though the analysis and
results are easily adapted to a general base $a$, simply by stipulating that
$p$ not divide $a$.

\subsection*{Patterns in Prime Divisors of Mersenne Numbers}

Let's look at patterns in divisors of Mersenne numbers, and see if Fermat's
broader theorem emerges as a conjecture. We will also see that these patterns
in divisors are what the third proposition is all about. Remember that
Fermat's theorem by itself is not obvious to observe; our contention is that
Fermat's broader theorem will be much easier to discover than just Fermat's
theorem on its own, and that the path of discovery flows in tandem with
discovering the third proposition. Note that we are going to look for patterns
in all prime divisors of all Mersenne numbers, not just in those Mersenne
numbers with prime exponents, nor only those with exponents directly relevant
to Fermat's theorem. This is a good example of how one doesn't always want to
restrict to as narrow a field of view as may be possible. Prime factorizations
through exponent $22$ are shown in Table \ref{TableKey}, all within Fermat's
easy grasp.%

\begin{table}[tbp] \centering
$%
\begin{array}
[c]{ll}%
M_{2}=3 & M_{13}=8191\\
M_{3}=7 & M_{14}=16\,383=3\times43\times127\\
M_{4}=15=3\times5 & M_{15}=32\,767=7\times31\times151\\
M_{5}=31 & M_{16}=65\,535=3\times5\times17\times257\\
M_{6}=63=3^{2}7 & M_{17}=131\,071\\
M_{7}=127 & M_{18}=262143=3^{3}7\times19\times73\\
M_{8}=255=3\times5\times17 & M_{19}=524\,287\\
M_{9}=511=7\times73 & M_{20}=1048\,575=3\times5^{2}11\times31\times41\\
M_{10}=1023=3\times11\times31 & M_{21}=2097\,151=7^{2}127\times337\\
M_{11}=2047=23\times89 & M_{22}=4194\,303=3\times23\times89\times683\\
M_{12}=4095=3^{2}5\times7\times13 & \,
\end{array}
$\caption{Prime factorizations of Mersenne numbers}\label{TableKey}%
\end{table}%

Now put Fermat's theorem and his propositions out of your mind, and let's just
look at the data through exponent $22$, as Fermat may have. We immediately see
some conjectures.

Once any number occurs as a divisor, it appears to propagate as a divisor of
precisely those Mersenne numbers whose exponents lie in a single simple
arithmetic progression. This provides the first proposition and the broader
claim of October, less Fermat's theorem itself.

\subsection*{The Divisibility Conjecture}

Another conjecture is that every odd prime $p$ occurs as a divisor, and that
this is certain to occur by a particular point in the sequence of Mersenne
numbers. Let's see where. Every prime divisor shown occurs for an exponent
less than the prime itself. And sometimes that happens just in the nick of
time, namely $p$ occurs first in $M_{p-1}$. Through our displayed range of
exponents, with $p-1\leq22$, this latter happens for $p=3,5,11,13,19$. The
other prime divisors through $p\leq23$, namely $7,17,23$, first occur earlier
than for exponent $p-1$. But even though they occur earlier than in the nick
of time, namely $p=7$ first in $M_{3}$, $p=17$ first in $M_{8}$, and $p=23$
first in $M_{11}$, we see that in all three cases they appear first for an
exponent $k$ that is a divisor of $p-1$ (namely $3\mid\left(  7-1\right)  $,
$8\mid\left(  17-1\right)  $, and $11\mid(23-1)$). Thus $p-1$ lies in the
simple arithmetic progression beginning with $k$, so it follows that $p$
divides $M_{p-1}$ in addition to dividing $M_{k}$. This leads to the
conjecture that $p$ always occurs in $M_{p-1}$ (Fermat's theorem), and
sometimes earlier, first in $M_{k}$ for some divisor $k$ of $p-1$. This is
exactly Fermat's October statement for $a=2$, and we'll call the claim there,
that the exponent $k$ (the first occurrence of $p$ as a divisor) must divide
$p-1$, the \emph{divisibility conjecture}, which implies Fermat's theorem (and
vice versa).

\begin{DC}
For any odd prime $p$, if $k$ is the least exponent for which $2^{k}-1$ is
divisible by $p$, then $k\mid\left(  p-1\right)  $. Moreover, such a $k$
always exists.
\end{DC}

Imagine how exciting these observations must have been for Fermat.

Even without further data, we (and Fermat) can easily go way beyond $p=23$ to
test the divisibility conjecture on all the prime factors visible in our
table, and thus confirm Fermat's theorem in more cases. For instance, the next
visible divisor is $p=31$, first occurring as $31\mid M_{5}$. Since the
divisibility conjecture $5\mid(31-1)$ holds in this case, it follows as well
that $31\mid M_{30}$. Thus, even though our table stops before $M_{30}$, we
have verified that Fermat's theorem is true for $p=31$. In general, for any
prime $p$ occurring in our data, we need merely ask, for its first occurrence
as a factor in some $M_{k}$, whether $k\mid(p-1)$, and this detects the truth
or falsity of Fermat's theorem for $p$.

\begin{exercise-unnumbered}
Confirm Fermat's theorem for $p=683$ and $524287$ just from the data above.
\end{exercise-unnumbered}

A distinction between two examples of this will illustrate an important point
for narrowing in on Fermat's third proposition and the amazing tool it
provided him for testing primality of Mersenne numbers. Consider the
occurrence of the primes $8191$ in $M_{13}$ and $683$ in $M_{22}$ in our data.
They should allow us to vet the truth (or falsity) of Fermat's theorem for
$p=8191$ and $683$, simply by verifying the divisibility conjecture. However,
to reliably state the divisibility conjecture for the occurrences at hand, we
should check first that $M_{13}$ and $M_{22}$ provide, respectively, the first
occurrences of these prime factors. Since each prime divisor occurs only for
exponents in a single simple arithmetic progression, the occurrence of $8191$
in $M_{13}$ is certainly the first, since the corresponding exponent $13$ is
prime. Therefore verifying $13\mid(8191-1)$ confirms the divisibility
conjecture, and thus the truth of Fermat's theorem for $p=8191$. However, in
the case of $683$, it could perhaps happen that its first occurrence is
earlier than in $M_{22}$, say in $M_{11}$, in which case the prime divisors of
$M_{22}$ might not definitively detect the veracity or falsity of their
instances of Fermat's theorem. Our table shows, though, that $M_{22}$ is the
first occurrence of $683$. So verifying $22\mid(683-1)$ confirms the
divisibility conjecture, and thus Fermat's theorem, for $p=683$. The point to
remember going forward is that if a prime occurs as a divisor of a Mersenne
number with a prime exponent, as in the occurrence of $8191$ in $M_{13}$, then
this is definitely its first occurrence, so the divisibility conjecture will
apply immediately to faithfully vet the truth or falsity of Fermat's theorem
for that prime.

\subsection*{The Form of Prime Divisors of Mersenne Numbers}

To summarize from the point just made, we see that if we look at a
(necessarily odd) prime divisor $p$ of a Mersenne number $M_{q}$ with odd
prime exponent $q$, and if we believe the divisibility conjecture (equivalent
to Fermat's theorem), then $q$ must divide $p-1$. In other words, $p$ must be
of the form $2jq+1$. Notice how this already came up as a possibility
discussed after Fermat's second proposition, but now it is confirmed, assuming
Fermat's broader theorem (which includes the divisibility conjecture). And lo
and behold, this is exactly what Fermat says in his third proposition!

\begin{FTP}
[June 1640]If $q$ is [an odd] prime, then the only possible prime divisors of
$M_{q}$ are numbers that are one greater than an even multiple of $q$.
\end{FTP}

\begin{center}%
\begin{center}
\includegraphics[
height=0.7422in,
width=3.8794in
]%
{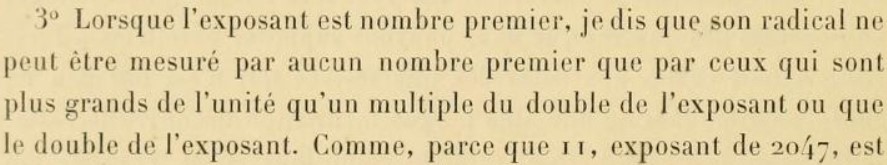}%
\\
Fermat's third proposition, June 1640 \cite[pp. 194--9]{fermatoeuvres}%
\end{center}

\end{center}

What a bombshell. Fermat has reduced the possible prime divisors of a Mersenne
number with prime exponent to a tiny fraction of all primes. Notice that to
obtain this we used Fermat's broader theorem, which combines Fermat's theorem
with the knowledge that a given prime $p$ occurs as a divisor of precisely the
Mersenne numbers in a simple arithmetic progression, whose beginning exponent
$q$ divides $p-1$, what we have called the divisibility conjecture. In fact
Fermat's third proposition is simply an alternative way of stating the
divisibility conjecture (and thus Fermat's theorem) in the case where the
beginning exponent $q$ is prime, which was Fermat's precise context in seeking
Mersenne primes.

It thus seems that Fermat knew his broader theorem for the base $2$ in June,
in particular the divisibility conjecture, and used it to obtain his third
proposition as a corollary targeted at divisors of Mersenne numbers with prime
exponent. How, indeed, could Fermat possibly have conjectured his third
proposition without appealing to his broader theorem?

Notice an easy generalization of Fermat's third proposition in the case where
$q$ is not necessarily prime.

\begin{GFTP}
If $p$ is a prime divisor of $M_{q}$ ($q$ not necessarily prime), but does not
occur as a divisor of any $M_{d}$ for $d$ a proper divisor of $q$, then $p$
must be of the form $jq+1$ ($2jq+1$ for $q$ odd).
\end{GFTP}

%

\begin{proof}%

By part b) of the Lemma above, the least exponent $k$ for which $p$ divides
$2^{k}-1$ must divide $q$. But our hypothesis prevents this $k$ from being a
proper divisor of $q$. Thus $k=q$, so by Fermat's broader theorem, $q$ must
divide $p-1$, hence $p$ has the claimed form.%

\end{proof}%

In a moment we'll find good use for this.

\subsection*{Advancing Toward $M_{37}$}

Let's use our methods to advance step by step beyond $M_{22}$ in the table
above, as could Fermat, to analyze $M_{n}$ for primality and divisors. We are
still aiming, with Fermat, towards $M_{37}$, in order to respond to Frenicle.

Fermat found that $47$ is a divisor of $M_{23}=8388607=\allowbreak
47\times178\,481$, so it is the second Mersenne prime imposter, after $M_{11}%
$. \cite[p. 210]{fermatoeuvres} His third proposition will produce this, but
this instance is also indicative of something more general.

\begin{exercise-unnumbered}
Show that if $p$ is a Sophie Germain prime (i.e., both $p$ and $q=2p+1$ are
prime), and if $p\equiv3\operatorname*{mod}4$, then $q$ is a divisor of
$M_{p}$. Hint: Is $2$ a square $\operatorname{mod}q$ (use quadratic reciprocity)?
\end{exercise-unnumbered}

Next, $M_{24}=16\,777\,215=3^{2}5\times7\times13\times17\times241$ is easily
obtained from our table for $M_{d}$ with $d\mid24$.

A nice illustration of our generalization above of Fermat's third proposition
to composite exponents is for $M_{25}=33\,554\,431$. We know from earlier only
that it is divisible by $M_{5}=31$, leaving a quotient $\allowbreak1082\,401$.
We can determine that the prime divisor $31$ has multiplicity only $1$ in
$M_{25}$, by checking that $31$ does not divide this quotient. The
generalization then tells us that all prime divisors of the quotient
$1082\,401$ must be of the form $50j+1$, making it not laborious to find by
division that $601$ is a prime divisor, leaving a final prime divisor $1801$.
So $M_{25}=33\,554\,431=31\times601\times1801.$

\begin{exercise-unnumbered}
Use Fermat's three propositions, our generalization of his third, and our
earlier table, to help you easily find, by hand, the prime factorizations of
$M_{28},M_{30},M_{36}.$ Also at least partly factor $M_{26},M_{27}%
,M_{29},M_{32},M_{33},M_{34},M_{35}.$
\end{exercise-unnumbered}

Here are the results of your exercise work, plus more information, especially
for $M_{29},M_{31}$.

$M_{26}=67\,108\,863=3\times2731\times8191.$

$M_{27}=134\,217\,727=7\times73\times262\,657.$

$M_{28}=268\,435\,455=3\times5\times29\times43\times113\times127$.

$M_{29}=536\,870\,911$ is a possible Mersenne prime. In 1738 Euler wrote that
$1103$ is a factor. \cite[v. 1, p. 17]{dickson}\cite[E26]{eulerarchive}%
\cite[p. 187]{SG-book} It is a little mysterious why he did not mention an
even smaller prime factor, namely $233$, which Fermat could easily have
discovered using his third proposition. The complete factorization is
$233\times1103\times2089.$ So $M_{29}$ is the third Mersenne prime imposter.

$M_{30}=1073\,741\,823=3^{2}7\times11\times31\times151\times331.$

$M_{31}=2147\,483\,647$ is another possible Mersenne prime. This is the number
that Frenicle jumped past in choosing the size of perfect number with which to
challenge Fermat, forcing Fermat on to $M_{37}$. So what happened to $M_{31}$?
In 1772 Euler wrote a letter to Bernoulli in which he claimed that Fermat had
proven that $M_{31}$ is prime, although I know of no evidence for this, and it
seems a stretch. Euler says that he proved it is prime himself, and that it is
the largest known prime at that time. \cite[v. 1, p. 19]{dickson}%
\cite[E461]{eulerarchive}\cite[p. 188]{SG-book} He had something up his sleeve
beyond Fermat's third proposition. He first says to Bernoulli that $M_{31}$
can have as prime divisors only numbers of the form $248n+1$ or $248n+63$, and
then says that he checked all these up to $46339$ (the square root of $M_{31}$
is $46340.\allowbreak9\dots$), none of which were found to be divisors, thus
confirming the primality of $M_{31}$. The corresponding perfect number is
$2^{30}M_{31}=2305\,843\,008\,139\,952\,128$, with $19$ digits, smaller than
Frenicle's challenge to Fermat.

So how did Euler obtain a twofold reduction in potential prime divisors by
restricting the possible forms from $62m+1$, according to Fermat's third
proposition, to $248n+1$ and $248n+63$? We'll think about this using more
modern concepts like primitive roots, congruence, and quadratic reciprocity;
Euler had his own ways of knowing all this. Euler knew that from $2^{31}-1$
being divisible by a prime $q$, he could conclude that $2$ is a square
$\operatorname{mod}q$ (Hint: Think primitive root.). And Euler knew that this
can only happen if $q\equiv\pm1\left(  \operatorname{mod}8\right)  $ (Hint:
Quadratic reciprocity). Thus $q=8l\pm1$. Combining this with $q=62m+1$ we have
that $q=248r+1$ or $248r+63$, as Euler claimed.

Imagine Euler checking the divisibility of $M_{31}$ by all the primes up to
$46339$ that are of these two residue forms (it seems Euler had a table of
primes up to $100,000$ \cite[E26]{eulerarchive})! He had to determine which of
the $4792$ primes less than $46339$ are in one of the two residue classes;
there are $84$ of them, beginning with the number $311$. And then he had to
test each of these as a divisor of $M_{31}$. Apparently all the results were negative.

\begin{exercise-unnumbered}
Estimate how long it took Euler to verify that $M_{31}$ is prime, even with
his twofold reduction of possible divisors.
\end{exercise-unnumbered}

Continuing,

$M_{32}=4294\,967\,295=3\times5\times17\times257\times65\,537.$

$M_{33}=8589\,934\,591=7\times23\times89\times599\,479.$

$M_{34}=17\,179\,869\,183=3\times43\,691\times131\,071.$

$M_{35}=34\,359\,738\,367=31\times71\times127\times122\,921.$

$M_{36}=68\,719\,476\,735=3^{3}5\times7\times13\times19\times37\times
73\times109.$

\subsection*{Fermat Responds to Frenicle's Challenge}

This brings our discussion finally to $M_{37}$, the crux of Frenicle's
challenge to Fermat, and a key stimulus for Fermat's discovery of his three
propositions and Fermat's theorem. So what does Fermat say about $M_{37}$? In
his June letter, immediately after divulging his three propositions, Fermat
writes that by using the shortcut from his third proposition, he can determine
that $M_{37}$ is divisible by $223$ \cite[v. 1, p. 12]{dickson}\cite[p.
199]{fermatoeuvres}\cite[ch. VI,\S II]{mahoney}\cite[ch. II,\S 4]{weil}.%

\begin{center}
\includegraphics[
height=1.637in,
width=3.9104in
]%
{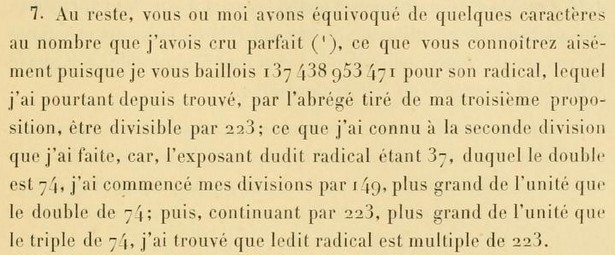}%
\\
Fermat's response to Frenicle's challenge, June 1640 \cite[p. 199]%
{fermatoeuvres}%
\end{center}

Looking at the primes one greater than an even multiple of $37$ as the only
possible divisors of $M_{37}$, he begins with $149=2\times74+1$, and continues
to $223=3\times74+1$, which he finds is a divisor. In fact the complementary
factor is also prime, although Fermat could not know this.

So $M_{37}=2^{37}-1$ $=137\,438\,953\,471=\allowbreak223\times616\,318\,177$,
and it does not lead to another perfect number. Voila, and bravo to Fermat!

\subsection*{What Ifs?}

Fermat's letters provide a rare and wonderful glimpse into the genesis of a
key mathematical discovery. There is just enough richness there to wonder also
how it might have turned out differently.

If, for example, Frenicle had challenged Fermat with $M_{31}$ rather than
$M_{37}$, Fermat would not have been able to confirm anything, even using his
third proposition, since $M_{31}$ is prime, but this is too large for Fermat
to have sussed. Would Fermat then even have bothered to tell Frenicle about
his second and third propositions? Or, what if $M_{37}$ had turned out to be
prime, or an imposter with no modestly sized factor? After all, the factor
$223$ is auspiciously small to encounter. Again Fermat would have little of
note to present to Frenicle, and might not have said anything. Thus is seems
that slightly different circumstances could have nixed Fermat's sharing of his
discoveries; and this might have affected his other important applications,
e.g., to the forms $a^{n}-1$ and $a^{n}+1$. Perhaps we can at least expect
that by the time Euler got involved many decades later, all would have been revealed.

\bigskip

\noindent\textbf{David Pengelley} is professor emeritus at New Mexico State
University, and courtesy professor at Oregon State University. His research is
in algebraic topology and history of mathematics. He develops the pedagogies
of teaching with student projects and with primary historical sources, and
created a graduate course on the role of history in teaching mathematics. He
relies on student reading, writing, and mathematical preparation before class
to enable active student work to replace lecture. He has received the MAA's
Deborah and Franklin Tepper Haimo teaching award, loves backpacking and
wilderness, is active on environmental issues, and has become a fanatical
badminton player.

\end{document}